\documentclass[12pt]{article}
\usepackage{amsmath,amssymb}
\newcommand{\lan}{\langle}
\newcommand{\ran}{\rangle}
\newcommand{\lt}{\left}
\newcommand{\rt}{\right}
\newcommand{\pf}{\mbox{\bf Proof}}
\newcommand{\ra}{\rightarrow}
\newcommand{\qed}{\hfill{\rule {.1in}{.1in}}}

\newcommand{\ol}{\overline}

\newcommand{\vp}{\varphi}
\newtheorem{thm}{Theorem}[section]

\newtheorem{prop}[thm]{Proposition}

\title{A criterion for finite rank $\lambda$-Toeplitz operators}
\author{Mark C. Ho\\
\scriptsize{Department of Applied Mathematics,
National Sun Yat-Sen University}\\
\scriptsize{Kaohsiung, Taiwan}\\
\scriptsize{Email: hom@math.nsysu.edu.tw}}
\date{\scriptsize{{\bf AMS classification}: 47A05, 47B38}}
\begin{document}
\maketitle
\begin{abstract}
Let $\lambda$ be a complex number in the closed unit disc $\ol{\Bbb
D}$, and $\cal H$ be a separable Hilbert space with the orthonormal
basis, say, ${\cal E}=\{e_n:n=0,1,2,\cdots\}$. A bounded operator
$T$ on $\cal H$ is called a {\em $\lambda$-Toeplitz operator} if
$\lan Te_{m+1},e_{n+1}\ran=\lambda\lan Te_m,e_n\ran$ (where
$\lan\cdot,\cdot\ran$ is the inner product on $\cal H$). The subject
arises naturally from a special case of the operator equation
\[
S^*AS=\lambda A+B,\ \mbox{where $S$ is a shift on $\cal H$},
\]
which plays an essential role in finding bounded matrix $(a_{ij})$
on $l^2(\Bbb Z)$ that solves the system of equations
\[\lt\{\begin{array}{lcc}
a_{2i,2j}&=&p_{ij}+aa_{ij}\vspace{.05in}\\
a_{2i,2j-1}&=&q_{ij}+ba_{ij}\vspace{.05in}\\
a_{2i-1,2j}&=&v_{ij}+ca_{ij}\vspace{.05in}\\
a_{2i-1,2j-1}&=&w_{ij}+da_{ij}
\end{array}\rt.
\]
for all $i,j\in\Bbb Z$, where $(p_{ij})$, $(q_{ij})$, $(v_{ij})$,
$(w_{ij})$ are bounded matrices on $l^2(\Bbb Z)$ and $a,b,c,d\in\Bbb
C$. It is also clear that the well-known Toeplitz operators are
precisely the solutions of $S^*AS=A$, when $S$ is the unilateral
shift. In this paper we verify some basic issues,
such as boundedness and compactness, for $\lambda$-Toeplitz
operators and, our main result is to give necessary and sufficient conditions for finite rank $\lambda$-Toeplitz operators.
\end{abstract}

\section{Introduction}

Let $\cal H$ be a separable Hilbert space with an orthonormal basis,
say, ${\cal E}=\{e_n:n=0,1,2,\cdots\}$. Given $\lambda\in\ol{\Bbb
D}=\{z\in\Bbb C:|z|\leq1\}$, a bounded operator $T$ is called a {\em
$\lambda$-Toeplitz operator} if $\lan
Te_{m+1},e_{n+1}\ran=\lambda\lan Te_m,e_n\ran$ (where
$\lan\cdot,\cdot\ran$ is the inner product on $\cal H$). In terms of
the basis $\cal E$, it is easy to see that the matrix representation
of $T$ is given by
\[\lt(\begin{array}{cccccc}
a_0&a_{-1}&a_{-2}&a_{-3}&a_{-4}&\cdots\\
a_1&\lambda a_0&\lambda a_{-1}&\lambda a_{-2}&\lambda a_{-3}&\ddots\\
a_2&\lambda a_1&\lambda^2a_0&\lambda^2a_{-1}&\lambda^2a_{-2}&\ddots\\
a_3&\lambda a_2&\lambda^2a_1&\lambda^3a_0&\lambda^3a_{-1}&\ddots\\
a_4&\lambda a_3&\lambda^2a_2&\lambda^3a_1&\lambda^4a_0&\ddots\\
\vdots&\ddots&\ddots&\ddots&\ddots&\ddots
\end{array}\rt)
\]
for some double sequence $\{a_n:n\in\Bbb Z\}$, and the boundedness
of $T$ clearly implies that $\sum|a_n|^2<\infty$. Therefore, it is
natural to introduce the notation
\[
T=T_{\lambda,\vp},
\]
where $\vp\sim\sum_{-\infty}^\infty a_ne^{in\theta}$ belongs to
$L^2=L^2(\Bbb T)$, the Hilbert space of square integrable functions
on the unit circle $\Bbb T$, with inner product
\[
\lan f,g\ran={1\over2\pi}\int_0^{2\pi}f{\ol g}d\theta,
\]
and consider $T_{\lambda,\vp}$ as an operator acting on $H^2\subseteq L^2$, the {\em Hardy space}
\[
\lt\{f\in L^2:\int_0^{2\pi}f(e^{i\theta})e^{in\theta}d\theta=0,\ n<0
\rt\}
\]
with the identification ${\cal H}=H^2$ and $e_n$ identified with the function $e^{in\theta}$, $n\geq0$. It is well-known that if $f\sim\sum_0^\infty a_ne^{i\theta}\in H^2$, then the analytic function
\[
\hat f(z):={1\over2\pi}\int_0^{2\pi}{f(e^{i\theta})\over1-ze^{-i\theta}}d\theta,\ |z|<1
\]
equals the analytic function defined by the series $\sum_0^\infty a_nz^n$, $|z|<1$, and, by a theorem of Fatou,
$\hat f_r(e^{i\theta})=\hat f(re^{i\theta})\ra f(e^{i\theta})$ a.e. $\theta$ and in $L^2$ (and for this reason $f$ is sometimes called the {\em boundary value} of $\hat f$). Hence $H^2$ is often identified with the Hilbert space of analytic functions $\{\hat f: f\in H^2\}$, with inner product
\[
\lan\hat f, \hat g\ran={1\over2\pi}\int_0^{2\pi}f\ol gd\theta,
\]
and we shall not make any distinction between $\hat f$ and $f$ throughout the rest of this paper. Also note that when $\lambda=1$ and $\vp\in L^\infty=L^\infty(\Bbb T)$, the matrix of $T_{1,\vp}$ is the matrix of the bounded Toeplitz operator $T_\vp f=P(\vp f)$, $f\in H^2$, where $P$ is the projection from $L^2$ on to $H^2$.
Here we refer the reader to \cite{RD:a} and \cite{KH:a}, both of which are excellent sources of information on the theory of Hardy spaces and Toeplitz operators on $H^2$.

The subject ``$\lambda$-Toeplitz operator" arises in a natural way
from the study of the bounded matrix $(a_{ij})$ on $l^2(\Bbb Z)$
that solves the system of equations
\[\lt\{\begin{array}{lcc}
a_{2i,2j}&=&p_{ij}+aa_{ij}\vspace{.05in}\\
a_{2i,2j-1}&=&q_{ij}+ba_{ij}\vspace{.05in}\\
a_{2i-1,2j}&=&v_{ij}+ca_{ij}\vspace{.05in}\\
a_{2i-1,2j-1}&=&w_{ij}+da_{ij}
\end{array}\rt.\eqno(\star)
\]
for all $i,j\in\Bbb Z$, where $(p_{ij})$, $(q_{ij})$, $(v_{ij})$,
$(w_{ij})$ are bounded matrices on $l^2(\Bbb Z)$ and $a,b,c,d\in\Bbb
C$ (See \cite{Ho:g}). One of the major steps for analyzing the solutions of ($\star$) is to determine bounded operator $A$ on $\cal H$ satisfying the operator equation below:
\[
S^*AS=\lambda A+B,
\]
where $S$ is a shift on $\cal H$, $B$ is fixed, and $|\lambda|\leq1$.
Notice that if we consider the map on ${\cal B}(\cal H)$:
\[
\phi(A)=S^*AS,\ A\in{\cal B}(\cal H),
\]
with $S$ being the unilateral shift, i.e., $Se_n=e_{n+1}$,
$n=0,1,2,\cdots$, then it is not difficult to see, by definition,
that the $\lambda$-Toeplitz operators are the ``eigenvectors" for
$\phi$ associated with $\lambda$. For instance, the Toeplitz
operators are just the ``eigenvectors" for $\phi$ associated with
the eigenvalue 1, and for every $\lambda$ with $|\lambda|=1$, the
``eigenvectors" of $\phi$ for the eigenvalue $\lambda$ are the {\em
Toeplitz-composition operators}, i.e., elements in the so-called
{\em Toeplitz-composition $C^*$-algebra}, which is a $C^*$-algebra generated by the Toeplitz algebra $\bf A$ and a single composition operator on $H^2$ (See, for example,
\cite{MJ:a} and \cite{KMM:a}) since
\[
T_{\lambda,\vp}=U_{\lambda}T_{\vp_{\ol\lambda,+}},
\]
where $T_{\vp_{\ol\lambda,+}}$ is the Toeplitz operator with symbol
\[
\vp_{\ol\lambda,+}\sim\sum_{n=-\infty}^\infty b_ne^{in\theta},\
b_n={\ol\lambda}^na_n\ \mbox{if}\ n\geq0\ \mbox{and}\ b_n=a_n\
\mbox{if}\ n<0,
\]
and the unitary
\[
U_\lambda e_n=\lambda^n e_n,\ n=0,1,2,\cdots
\]
is clearly the same with the composition operator on $H^2$ induced
by $\lambda z$ (For information on the theory of composition operators on Hardy spaces, we refer the readers to \cite{CM:a})). On the other hand, when $|\lambda|<1$, $\lambda$-Toeplitz operators are no longer Toeplitz-composition operators, but they can be written as the sum of weighted composition operators and their adjoints (See the discussion below following Proposition \ref{bound2}).

In this paper we study some basic properties, such as boundedness and compactness, for $\lambda$-Toeplitz
operators and while doing so, we improve a result concerning the compactness for certain weighted composition operators (Proposition \ref{weighted}). Our main result, on the other hand, is to give necessary and sufficient conditions for the parameter $\lambda$ and the symbol $\vp$ so that the corresponding $\lambda$-Toeplitz operator has finite rank (Proposition \ref{finite}).
\section{Bounded and compact $\lambda$-Toeplitz operators}
\label{compact}
We shall proceed in two cases: $|\lambda|=1$ and $|\lambda|<1$.\vspace{.1in}\\
$|\lambda|=1$:\vspace{.05in}

For $\vp\sim\sum a_ne^{in\theta}$ in $L^2$, let us recall the unitary $U_\lambda$, the function
\[
\vp_{\ol\lambda,+}\sim\sum_{n=-\infty}^\infty b_ne^{in\theta},\
b_n={\ol\lambda}^na_n\ \mbox{if}\ n\geq0\ \mbox{and}\ b_n=a_n\
\mbox{if}\ n<0,
\]
defined previously and the relation
$T_{\lambda,\vp}=U_{\lambda}T_{\vp_{\ol\lambda,+}}$.
Hence by the classical results of Toeplitz operators, we have
\begin{prop}
Let $|\lambda|=1$ and $\vp\sim\sum a_ne^{in\theta}$ in $L^2$.
\begin{itemize}
\item[1.] The
the $\lambda$-Toeplitz operator $T_{\lambda,\vp}$ is bounded if and
only if the Toeplitz operator $T_{\vp_{\ol\lambda,+}}$ is bounded
or, equivalently, $\vp_{\ol\lambda,+}\in L^\infty=L^\infty(\Bbb T)$.
Moreover, in this case, we have
\[
\|T_{\lambda,\vp}\|=\|\vp_{\ol\lambda,+}\|_\infty.
\]
\item[2.] $T_{\lambda,\vp}$ is compact if and only if $\vp_{\ol\lambda,+}\equiv0$, i.e., if and only if $\vp\equiv0$.
\end{itemize}
\label{bound1}
\end{prop}
{\bf Remark.} We also have
$\|T_{\lambda,\vp}\|=\|T_{\lambda,\vp}U_{\ol\lambda}\|=\|T_{\vp_{\ol\lambda,-}}\|$
(hence $\|T_{\lambda,\vp}\|=\|\vp_{\ol\lambda,-}\|_\infty$), where
$T_{\vp_{\ol\lambda,-}}$ is the Toeplitz operator with symbol
\[
\vp_{\ol\lambda,-}\sim\sum_{n=-\infty}^\infty b_ne^{in\theta},\
b_n=a_n\ \mbox{if}\ n\geq0\ \mbox{and}\ b_n={\ol\lambda}^na_n\
\mbox{if}\ n<0.
\]
Note also that the boundedness of $\vp$ alone, not like the case in
the Toepliz operator $T_\vp$, does not guarantee the boundedness of
$T_{\lambda,\vp}$. Here is an example: Consider the $2\pi$-periodic
extension of $\vp(\theta)=\theta$, $0\leq\theta<2\pi$ with Fourier
series
\[
1+\sum_{n\not=0}{1\over in}e^{in\theta}.
\]
Obviously $\vp\in L^\infty$. But it is not difficult to see that
$T_{-1,\vp}$ is not bounded since the function $\vp_{-1,+}\not\in
L^\infty$. In fact in this case, we have $T_{\lambda,\vp}$ is
bounded (or, equivalently, $\vp_{\ol\lambda,+}\in L^\infty$) if and
only if $\lambda=1$.\vspace{.1in}\\
$|\lambda|<1$:\vspace{.05in}

The boundenes and compactness of
$T_{\lambda,\vp}$ in this case are, basically, automatic:
\begin{prop}
Let $|\lambda|<1$ and $\vp$ is measurable on $\Bbb T$. Then the
$\lambda$-Toeplitz operator $T_{\lambda,\vp}$ is bounded if and only
if $\vp\in L^2$. In fact, every bounded $\lambda$-Toeplitz operator
is Hilbert-Schmidt if $|\lambda|<1$, and we have
\[
\|T_{\lambda,\vp}\|\leq(1-|\lambda|^2)^{-1/2}(\sum_{n=-\infty}^\infty|a_n|^2)^{1/2},
\]
where $\vp\sim\sum a_ne^{in\theta}$. \label{bound2}
\end{prop}
\pf\ Just observe that when $|\lambda|<1$,
\[
(1-|\lambda|^2)^{-1}(\sum_{n=-\infty}^\infty|a_n|^2)<\infty
\]
if and only if $\vp\in L^2$.\qed\vspace{.05in}

We will see later that $T_{\lambda,\vp}$ is actually in the trace class if $|\lambda|<1$. The reason for the redundancy here is that we would like to compare what we have so far to a similar result from another class of operators on $H^2$.
Let $\tau$ be an analytic map from $\Bbb D$ into $\Bbb D$. Then,
given $\psi$ analytic on $\Bbb D$, the weighted composition operator
$W_{\psi,\tau}$ on $H^2$ is defined by
\[
W_{\psi,\tau}f:=\psi\cdot(f\circ\tau).
\]
Weighted composition operators on Hardy spaces have been
receiving increasing attention from operator theorists in the past
decade on subjects including boundedness, compactness, spectrum
(see, for example, \cite{CA:a}, \cite{CA:b}, \cite{CE:a},
\cite{GG:a}, \cite{GG:b}, \cite{MZ:a}, \cite{SH:a}). In particular,
on the subject of boundedness and compactness, results so far have
suggested that $W_{\psi,\tau}$ may be bounded, or even compact with
unbounded $\psi$. For example, G. Gunatillake proves that
$W_{\psi,\tau}$ is bounded and in fact, compact, if $\ol{\tau(\Bbb D)}\subseteq\Bbb D$ ($\ol{\tau(\Bbb D)}$ being the closure of $\tau(\Bbb D)$) and
$\psi\in H^2$ by showing $\|W_{\psi,\tau}|_{z^nH^2}\|$ tends to zero
as $n\ra\infty$ (see Theorem 2, \cite{GG:b}). This of course implies
immediately that $T_{\lambda,\vp}$ is compact if $|\lambda|<1$ since both $W_{\vp_+,\lambda z}$ and $W_{\ol\vp_-,\ol\lambda z}$ are compact and
\[
T_{\lambda,\vp}=W_{\vp_+,\lambda z}+W_{\ol\vp_-,\ol\lambda z}^*,
\]
where $\vp_+=P\vp$ and $\ol\vp_-$ is the $H^2$ function whose boundary value is given by the conjugate of $\vp_-=(I-P)\vp$, i.e.,
\[
\ol\vp_-(z)=\sum_{n=1}^\infty\ol a_{-n}z^n,\ |z|<1.
\]
So it would appear that $T_{\lambda,\vp}$ being compact if $|\lambda|<1$ is just a special case and therefore a predictable result of Gunatillake's work, and its being Hilbert-Schmidts is merely a coincidence under special circumstances.

However, we will show that these ``special circumstances" are not as special as one would think. In fact, we will show that the ``Hilbert-Schmidt" conclusion also holds under Gunatillake's more general assumptions:
\begin{prop}
Suppose that $\ol{\tau(\Bbb D)}\subseteq\Bbb D$. Then
$W_{\psi,\tau}$ is Hilbert-Schmidt if $\psi\in H^2$.
\label{weighted}
\end{prop}
\pf\ The key is to express $W_{\psi,\tau}$ in the form of an integral operator and analyze its kernel.
Let us consider the reproducing kernel
$K_\alpha(z)=(1-{\ol\alpha}z)^{-1}$ on $H^2$. Then, by the
assumption that $\ol{\tau(\Bbb D)}\subseteq\Bbb D$, the function
\[
\psi(e^{i\phi})\ol{K_{\tau(e^{i\phi})}(e^{i\theta})}=
\psi(e^{i\phi})(1-e^{-i\theta}{\tau(e^{i\phi})})^{-1}\in L^2(\Bbb
T\times\Bbb T).
\]
Then since for $f\in H^2$,
\[ (W_{\psi,\tau}f)(z)=\psi(z)f(\tau(z))=\psi(z)\lan
f,K_{\tau(z)}\ran,\ |z|<1,
\]
it is not difficult to recognize that $W_{\psi,\tau}$ is actually
the same with the restriction of the $L^2$ integral operator
\[
f(e^{i\phi})\longrightarrow{1\over2\pi}\int_0^{2\pi}
\psi(e^{i\phi})\ol{K_{\tau(e^{i\phi})}(e^{i\theta})}f(e^{i\theta})d\theta
\]
on the boundary values of functions in $H^2$. This completes the
proof since it is well-known that integral operators on $L^2(\Bbb
T)$ with square integrable kernels are Hilbert-Schmidts.\qed\vspace{.05in}\\
{\bf Remark.} With the above conventions, $T_{\lambda,\vp}$ can be
regarded as the integral operator
\[
f(e^{i\phi})\longrightarrow{1\over2\pi}\int_0^{2\pi}
(\vp_+(e^{i\phi})+\vp_-(e^{i\theta}))\ol{K_{\lambda
e^{i\phi}}(e^{i\theta})}f(e^{i\theta})d\theta,\ f\in H^2.
\]

We now finish this section with
\begin{thm}
$T_{\lambda,\vp}$ belongs to the trace class if $|\lambda|<1$.
\label{trace}
\end{thm}
\pf\ Let us consider the block form for the matrix of $T_{\lambda,\vp}$:
\[\lt(\begin{array}{cc}
A_m&B_m\\
C_m&D_m
\end{array}\rt),
\]
where $A_m$ is the upper-left $m\times m$ minor of the matrix.
Since $\|D_m\|=|\lambda|^m\|T_{\lambda,\vp}\|$ and the rank of
\[\lt(\begin{array}{cc}
A_m&B_m\\
C_m&0
\end{array}\rt)
\]
is at most $2m$, we see that the ($2m+1$)-th singular value of
$T_{\lambda,\vp}$ is less than or equal to
$|\lambda|^m\|T_{\lambda,\vp}\|$. This completes the
proof since the sequence of the sigular values is decresing.\qed
\section{$\lambda$-Toeplitz operators with finite rank}
Basically the results in Section \ref{compact} tells us that if $|\lambda|=1$, then $T_{\lambda,\vp}$ is never compact except in the trivial case and if, on the other hand, $|\lambda|<1$, then any choice of $\vp$ in $H^2$ will make $T_{\lambda,\vp}$ compact. Let us now consider a question in a some what similar nature: Given $|\lambda|<1$, how can we tell whether $T_{\lambda,\vp}$ is of finite rank or not by the choice of $\vp$? We answer this question in the following:
\begin{prop}
Let $|\lambda|<1$ and $\vp\sim\sum a_ne^{in\theta}\in L^2$.
Then $T_{\lambda,\vp}$ is of finite rank if and only if $\vp=0$ or
$\lambda=0$. \label{finite}
\end{prop}
\pf\ Since $T_{0,\vp}$ has rank at most two, we will assume, from now on, that $0<|\lambda|<1$.

First we show that if $\vp_+=0$ or $\vp_-=0$, i.e., if the matrix of $T_{\lambda,\vp}$ is either lower triangular or upper triangular, then $\dim{\cal R}(T_{\lambda,\vp})=\infty$ unless $\vp=0$ (here ${\cal R}(A)$ is the range of $A$). Let us assume, without loss of generality, $\vp\not=0$ and $\vp_-=0$. This means that $\vp_+\not=0$ and therefore
$T_{\lambda,\vp}$ is the weighted composition operator $W_{\vp_+,\lambda z}$, and hence $S^{*n_0}T_{\lambda,\vp}$ is the weighted composition operator $W_{\psi,\lambda z}$, where $n_0=\min\{n:a_n\not=0\}$ and $\psi=S^{*n_0}\vp_+$. Now suppose that $\dim{\cal R}(T_{\lambda,\vp})<\infty$. Then $\dim{\cal R}(W_{\psi,\lambda z})<\infty$. This leads to a contradiction since the spectrum of $W_{\psi,\lambda z}$, according to Theorem 1 in \cite{GG:a}, is the set
\[
\{\lambda^k\psi(0):k=0,1,2,\cdots\}\cup\{0\}
\]
and this set is infinite since $\psi(0)=a_{n_0}\not=0$ and $\lambda\not=0$.

What follows from the above discussion is that if $a_n\not=0$ for only finitely many nonnegative or negative $n$, then $\dim{\cal R}(T_{\lambda,\vp})=\infty$, unless $a_n=0$ for all $n$ (i.e. $\vp=0$). Indeed, for example, if $a_n=0$ for all but finitely many positive $n$ and $n_0=\max\{n>0:a_n\not=0\}\geq0$, then the matrix of the $\lambda$-Toeplitz operator $S^{*n_0}T_{\lambda,\vp}$ is upper triangular with nonzero diagonal, and hence $\dim{\cal R}(T_{\lambda,\vp})=\infty$. Consequently, $T_{\lambda,\vp}$ can not be of finite rank if at least one of $\vp_+$ and $\vp_-$ is a nonzero trigonometric polynomial.

Therefore, we shall proceed by assuming that $\vp_+$ and $\vp_-$ are nonzero and none of them are trigonometric polynomials. Since $T_{\lambda,\vp}$ has finite rank, the kernel of $T_{\lambda,\vp}$ contains nonzero polynomials (in $H^2$). Hence by the identity
\[
T_{\lambda,\vp}=W_{\vp_+,\lambda z}+W_{\ol\vp_-,\ol\lambda z}^*,
\]
we can find polynomials $p,\ q\not=0$ such that
\[
\vp_+(z)p(\lambda z)=q(\ol\lambda z),\ |z|<1\ \ (q=-P({\vp}_-p)).
\]
Notice the above equality also depends on the fact that the composition operator $C_{\ol\lambda z}$ and $P$ commute. It follows that $\vp_+$ is a
rational function. However, since $\vp_+$ is also in $H^2$, all poles of $\vp_+$ must lie strictly outside the unit disc (i.e., in the region $\{z:|z|>1\}$). So $\vp_+$ is actually continuous on the closed disc $\ol{\Bbb D}$. By considering $T_{\lambda,\vp}^*$, wee see that the same also holds for $\ol\vp_-$.

Now, without loss of generality, we assume that $\lambda>0$ since
otherwise, if $\lambda=|\lambda|e^{i\alpha}$, we have
$T_{\lambda,\vp}=U_{e^{i\alpha}}T_{|\lambda|,\vp_{e^{-i\alpha,+}}}$, and note that
\[
\vp_+(z)=(\vp_{e^{-i\alpha,+}})_+(e^{i\alpha}z).
\]
Let $\tilde\vp$ be the function on $\Bbb T$ defined by
\[
\tilde\vp(e^{i\theta})=\vp_+(e^{i\theta})+\ol{\ol\vp_-(\lambda e^{i\theta})},\ e^{i\theta}\in\Bbb T.
\]
So $\tilde\vp$ is continuous on $\Bbb T$, and hence $T_{\tilde\vp}$ is a bounded Toeplitz operator. On the other hand, again by the fact that $C_{\lambda z}$ commute with $P$, it is not difficult to check that (remembering $\lambda>0$)
\[
T_{\lambda,\vp}=T_{\tilde\vp}C_{\lambda z}.
\]
But because $C_{\lambda z}$ is one to one, the fact that $T_{\lambda,\vp}$ has finite rank implies that $T_{\tilde\vp}$ also has finite rank, hence $\tilde\vp=0$ and, consequently, $\vp=0$.
\qed

Mark C. Ho\\
Department of Applied Mathematics\\
National Sun Yat-Sen University\\
Kaohsiung, Taiwan 80424\\
hom@math.nsysu.edu.tw

\end{document}